\documentclass[a4paper]{article}
\newtheorem{thm}{Theorem}[section]
\newtheorem{lem}[thm]{Lemma}
\newtheorem{eg}[thm]{Example}
\newtheorem{prop}[thm]{Proposition}
\newtheorem{cor}[thm]{Corollary}
\newtheorem{defn}[thm]{Definition}
\newtheorem{rem}[thm]{Remark}

\newtheorem{rem-eg}[thm]{Remark and Example}

\newcommand{\smnoind}{{\smallskip\noindent}}
\newcommand{\id}{{\rm id}}
\newcommand{\Rep}{{\rm Rep}}

\begin{document}
\title{Duality of Hopf $C^*$-algebras}
\author{Chi-Keung Ng}
\date{}
\maketitle
\begin{abstract}
In this paper, we study the duality theory of Hopf $C^*$-algebras in a general 
``Hilbert-space-free'' framework. 
Our particular interests are the ``full duality'' and the ``reduced duality''. 
In order to study the reduced duality, we define the interesting notion of Fourier algebra of a general Hopf $C^*$-algebra.
This study of reduced duality and Fourier algebra is found to be useful in the study of other aspects of Hopf $C^*$-algebras (see e.g. \cite{Ng3}, \cite{Ng4} and \cite{Ng5}). 
\end{abstract}
\noindent{\small 2000 AMS Mathematics Classification numbers: Primary 46L55; 46L89; Secondary 22D15}
\par\medskip
\par\medskip
\par\medskip

\section*{0\hspace{0.25in}Introduction}
\par\bigskip

The study of Hopf $C^*$-algebras began with group 
$C^*$-algebras.
There is a strong relation between locally compact groups and 
commutative (or co-commutative) Hopf $C^*$-algebras (see \cite{Ior} and \cite{Val}).
The general (i.e$.\!$ non-commutative and non-co-commutative) Hopf $C^*$-algebras are getting more attention with the different studies of ``locally compact quantum groups'' in the recent literatures (see e.g. \cite{Baaj}, \cite{BS}, \cite{BMT}, \cite{DQV}, \cite{ES}, \cite{KV1}, \cite{KV}, \cite{VaV}, \cite{VV} and \cite{Wor3}). 

\bigskip

One important and interesting feature of Hopf $C^*$-algebras is their duality theory (which is also related to the interesting notion of ``amenability'' -- see the last paragraph of this introduction). 
In the literatures, duality is usually tied up with the notion of \emph{multiplicative unitaries}. 
The aim of this paper is to study the duality theory of Hopf $C^*$-algebras without referring explicitly to multiplicative unitary. 
The advantage of this approach is that there may not be a multiplicative unitary in the first place and different multiplicative unitaries may give ``the same duality'' between two Hopf $C^*$-algebras. 
Moreover, it is not clear the kind of conditions (e.g. regularity in \cite{BS} and manageability in \cite{Wor3}) required on those multiplicative unitaries (e.g. the tensor product of a regular multiplicative unitary with a manageable one may not be regular nor manageable any more). 

\bigskip

In this article, we will give a general intrinsic framework for the duality of Hopf $C^*$-algebras (which includes all sort of interesting multiplicative unitaries) and will pay more attention to the ``full duality'' (or ``universal duality'') and the ``reduced duality''. 
The first one will be studied in Section 2 (where its existence and related topics will be discussed). 
To define the second duality, we will needed some preparation. 
In particular, we need to have some understanding of the general duality framework, namely, the Fourier duality.  
Section 3 is devoted to the study of this. 

\bigskip

The reduced duality will be studied in Section 4. 
In order to define this, a natural notion of the ``Fourier algebra'' of a Hopf $C^*$-algebra is required.
This object is interesting as its definition is intrinsic (i.e. concerning just the structure of the given Hopf $C^*$-algebra without referring to any extra structure) but it is the ``core'' of the duality theory given by multiplicative unitaries. 
Using this, we obtained a kind of uniqueness result. 

\bigskip

The study of duality theory in this paper is important for the understanding of ``amenability'' of general Hopf $C^*$-algebras. 
In fact, roughly speaking, we can think of amenable Hopf $C^*$-algebras as those whose ``full duality'' and ``reduced duality'' coincide. 
We gave a detail study of this in another paper (\cite{Ng3}).
Moreover, the Fourier algebras and its relation with reduced duality (as given in Section 4) is also important for the study of amenability of Kac algebras and Kac systems (see \cite{Ng4} and \cite{Ng5}).
\par\medskip
\par\medskip
\par\medskip

\section{Preliminary and notations}
\par\bigskip

We begin with the definition of Hopf $C^*$-algebras and Hopf von Neumann algebras. 

\par\medskip

\begin{defn}\label{1.6} 
Let $S$ be a $C^*$-algebra.
\par\smnoind
(a) A non-degenerate $*$-homomorphism from $S$ to $M(S\otimes S)$ 
is said to be a \emph{comultiplication} if $(\delta\otimes \id)\delta= 
(\id\otimes\delta)\delta$. 
A $*$-homomorphism $\varphi\in S^*\setminus \{0\}$ is called a \emph{coidentity} for $\delta$ if $(\varphi \otimes \id)\delta = \id = (\id\otimes \varphi)\delta$. 
\par
\smnoind
(b) $S$ is called a \emph{pre-Hopf $C^*$-algebra} if it has a
comultiplication $\delta$. 
Moreover, $(S,\delta)$ is said to be a \emph{Hopf $C^*$-algebra}, if $\delta(S)(1\otimes S)\subseteq S\otimes S$ and $\delta(S)(S\otimes 1)\subseteq S\otimes S$.
\par
\smnoind
(c) Let $A$ be a $C^*$-algebra and $S$ be a pre-Hopf $C^*$-algebra. 
Then a non-degenerate $*$-homomorphism $\epsilon$ from $A$ to 
$M(A\otimes S)$ is called a \emph{coaction} if 
$(\epsilon\otimes \id)\circ\epsilon=(\id\otimes \delta)\circ\epsilon$ and $\epsilon(A)(1\otimes S)\subseteq A\otimes S$. 
\end{defn}
\par\medskip

Similarly, a non-degenerate $*$-homomorphism from $S$ to 
$M(S\otimes_{\max} S)$ is called a \emph{full comultiplication} if it satisfies the equality in Definition \ref{1.6}(a). 
$S$ is said to be a \emph{full Hopf $C^*$-algebra} if it has a comultiplication
satisfying similar conditions as in Definition \ref{1.6}(b). 

\bigskip

\begin{rem}\label{1.3}
(a) Similarly, we can define \emph{Hopf von Neumann algebras} and their \emph{coactions}. 

\smnoind
(b) Let $A$ and $B$ be $C^*$-algebras and let $\cal A$ and $\cal B$ 
be von Neumann algebras. 
Suppose that $\varphi$ and $\psi$ are $*$-homomorphisms from $A$ 
and $B$ to $\cal A$ and $\cal B$ respectively such that they are 
\emph{non-degenerate} (in the sense that the image of any approximate unit converges to $1$ weakly). 
Then $\varphi\otimes \psi$ extended to a map on $M(A\otimes B)$
Note that for any $F\in ({\cal A}\overline\otimes {\cal B})_*$ and 
$m\in M(A\otimes B)$, we have $(\varphi\otimes\psi)^*(F\cdot 
(\varphi\otimes\psi)(m)) = (\varphi\otimes\psi)^*(F)\cdot m$.

\smnoind
(c) By part (b) above, we see that 
if $A$ is a $C^*$-algebra and $S$ is a Hopf $C^*$-algebra, then
$S^{**}$ is a Hopf von Neumann algebra and 
any coaction $\epsilon$ on $A$ by $S$ induces a weakly
continuous coaction $\epsilon^{**}$ on $A^{**}$ by $S^{**}$.
\end{rem}

\par\medskip

For any subspace $N\subseteq S^*$, we denote $N^{\top} = \{x\in S: f(x) = 0$ for all $f\in N \}$. 

\par\medskip

\begin{prop}\label{1.7}
Let $(S,\delta)$ be a Hopf $C^*$-algebra 
and $e$ be a non-trivial central projection of $S^{**}$. 
Let $T = S/(e\cdot S^*)^\top$ with the canonical quotient map $Q$. 
Then the following statements are equivalent.

\smnoind
(1) $\ker (Q)\subseteq \ker (Q\otimes \id)\delta$;

\smnoind
(2) $\delta$ induces a coaction $\epsilon$ on $T$ by $S$ (i.e. 
$\epsilon\circ Q = (Q\otimes \id)\circ\delta$);

\smnoind
(3) $\overline{e\cdot S^*}^{\sigma(S^*,S)}$ is a right ideal of $S^*$.
\end{prop}
\par
\noindent {\bf Proof:}
It is clear that $\ker (Q)\subseteq \ker (Q\otimes \id)\circ\delta$ 
implies that $\delta$ induces a coaction on $T$ by $S$. 
If $\delta$ induces a coaction $\epsilon$ on $T$ by $S$, then for any 
$f\in T^*$ and $g\in S^*$, $Q^*(f)g = (f\otimes g)(Q\otimes \id)
\circ\delta = Q^*(h)$ where $h=(f\otimes g)\circ\epsilon$. 
Hence $\overline{e\cdot S^*}^{\sigma(S^*,S)}=Q^*(T^*)$ is a right ideal of $S^*$.
Finally, assume that $\overline{e\cdot S^*}^{\sigma(S^*,S)}$ is a right ideal
of $S^*$.
For any $x\in \ker (Q)$, we have $(Q^*(f)g)(x) = 0$ (for $f\in T^*$ and $g\in S^*$) and thus, $(Q\otimes \id)\delta(x) =0$. 
\par\medskip

\begin{rem}\label{1.8}
(a) There is a corresponding result for 
$\ker (Q)\subseteq \ker (\id\otimes Q)\delta$.
\par\smnoind
(b) By a similar argument as that of Proposition \ref{1.7}, we can 
also show that the followings are equivalent: 
\par\hspace{-1em}
(1) $\ker (Q)\subseteq \ker (Q\otimes Q)\delta$;
\par\hspace{-1em}
(2) $\delta$ induces a quotient Hopf $C^*$-algebra structure on $T$;
\par\hspace{-1em}
(3) $\overline{e\cdot S^*}^{\sigma(S^*,S)}$ is a subalgebra of $S^*$.
\par\smnoind
(c) For any locally compact group $G$, the only quotient Hopf $C^*$-algebra of 
$C^*_r(G)$ that has a compatible coaction by $C^*_r(G)$ is $C^*_r(G)$ 
itself (combine part (b) and Proposition \ref{1.7} with \cite[2.4]{DR}). 
\end{rem}

\par\medskip

\begin{prop}\label{1.9}
Let $\cal M$ be a Hopf von Neumann algebra with a comultiplication 
$\bar\delta$ and $e\in \cal M$ be a non-trivial central projection. 
Let $\bar Q$ be the canonical map from $\cal M$ to $e\cal M$. 
Then the following statements are equivalent.

\smnoind
(1) $\ker ({\bar Q})\subseteq \ker ({\bar Q}\otimes \id)
\bar \delta$;

\smnoind
(2) $(\bar Q\otimes \id)\circ\bar \delta\mid_{e{\cal M}}$ is 
a coaction on $e\cal M$ by $\cal M$;

\smnoind
(3) $\bar \delta(e)\geq e\otimes 1$; 

\smnoind
(4) $e{\cal M}_*$ is a right ideal of ${\cal M}_*$.
\end{prop}
\par
\noindent {\bf Proof:}
Assume that (1) holds. 
Then there exists a weak-*-continuous map $\bar {\epsilon}$
such that $\bar {\epsilon}\circ{\bar Q} = 
({\bar Q}\otimes \id)\circ\bar \delta$. 
It is clear that $\bar {\epsilon}(e) = \bar\epsilon (\bar Q(1)) = e\otimes 1$ and 
$\bar {\epsilon} = (\bar {Q}\otimes \id)\circ\bar \delta
\mid_{e{\cal M}}$ is a coaction. 
Hence (2) is established. 
It is easy to see that (2) implies (3) (as $({\bar Q}\otimes \id)\circ\bar \delta (e) = e\otimes 1$). 
Suppose that (3) is true (i.e. $(e\otimes 1)\bar {\delta}(1-e) = 0$). 
For any $z\in {\cal M}$ and $\omega, \nu\in {\cal M}_*$, we have 
$(1-e)(e\omega \cdot\nu)(z) = (\omega\otimes\nu)(e\otimes 1)
\bar {\delta}((1-e)z) = 0$ and thus, (4) holds. 
Finally, we need to show that (4) implies (1). 
Let $z\in \ker(\bar {Q})$ (i.e. $ez=0$).
For any $\omega\in e{\cal M}_*$ and $\nu\in {\cal M}_*$, one has 
$(\omega\otimes \nu)(\bar {Q}\otimes \id)\bar {\delta}(z) = 
(e\omega\cdot\nu)(z) = 0$. 
This completes the proof. 
\par\medskip

\begin{rem}\label{1.10}
(a) The equivalence of (3) and (4) is, in fact, \cite[2.2]{DR} in the case of 
central projections.
\par\smnoind
(b) We also have the von Neumann algebra version of Remark \ref{1.8}(b). 
In particular, the map $\bar \delta_e$ from $e{\cal M}$ to 
$e{\cal M}\otimes e{\cal M}$ defined by $\bar \delta_e(x) = 
(e\otimes e)\delta (x)$ is a comultiplication on $e\cal M$ if 
and only if $\bar \delta (e) \geq e\otimes e$. 
It is the case if and only if $e{\cal M}_*$ is a subalgebra of 
${\cal M}_*$. 
\par\smnoind
(c) Let $T$ be a quotient $C^*$-algebra of a Hopf $C^*$-algebra $S$ and $Q$ be
the canonical quotient map. 
Then $Q^*(T^*)$ is both $S$-invariant and $\sigma(S^*,S)$-closed. 
If $e_T$ is a central projection in $S^{**}$ such that $Q^*(T^*)=e_T\cdot S^*$ (as $Q^*(T^*)$ is $S$-invariant), then Propositions \ref{1.7} and \ref{1.9} as well as Remark \ref{1.8}(b) and part (b) above give relations between the properties of $T$ and the properties of
$e_T$. 
\end{rem}
\par\medskip

We recall the following notations from \cite{Wor2}: 
for any Hilbert space $H$ and any $C^*$-algebra $A$, let $C^*(H)$ be the 
collection of all non-degenerate separable $C^*$-subalgebras of ${\cal L}(H)$ 
and $\Rep(A,H)$ be the collection of all non-degenerate representations 
of $A$ on $H$. 

\par\medskip

\begin{defn}\label{1.11} (\cite[4.1]{Wor2}) 
Let $S$ and $T$ be two separable $C^*$-algebras and $X\eta (T\otimes S)$ (please refer to \cite{Wor2} for its meaning). 
Then $S$ is said to be \emph{generated by $X$} if for any separable Hilbert 
space $H$, any $B\in C^*(H)$ and any $\pi\in \Rep(S,H)$, condition (i)
below implies condition (ii):

\smnoind
(i) $(\id\otimes \pi)(X)\eta (T\otimes B)$;

\smnoind
(ii) $\pi(S)B$ is dense in $B$. 
\end{defn}
\par\medskip

The following simple result can be found in the second paragraph after \cite[4.1]{Wor2}.
\par\medskip

\begin{lem}\label{1.13} 
Let $H$ be a separable Hilbert space and $C$ be a separable $C^*$-algebra. 
If $X\eta (C\otimes {\cal K}(H))$, then there exists at most one $A\in 
C^*(H)$ such that $X\eta (C\otimes A)$ and $A$ is generated by $X$.
\end{lem}
\par\medskip

\begin{prop}\label{1.14}
Let $A$, $B$, $C$ and $D$ be separable $C^*$-algebras. 
Let $X\eta (C\otimes A)$ such that $A$ is generated by $X$.
\par\smnoind
(a) Suppose that there is a surjective $*$-homomorphism $\varphi$ 
from $A$ to $D$.
Then $D$ is generated by $(\id\otimes \varphi)(X)$. 
\par\smnoind
(b) Suppose that $\phi$ is a non-degenerate $*$-homomorphism from $A$ to
$M(B)$ such that $B$ is generated by $(\id\otimes \phi)(X)$. 
Then $\phi(A) = B$. 
\end{prop}
\par
\noindent {\bf Proof:}
(a) Let $H$ be a Hilbert space. 
Suppose that $E\in C^*(H)$ and $\pi\in \Rep(D,H)$ such that 
$(\id\otimes \pi)(\id\otimes \varphi)(X)\eta (C\otimes E)$. 
Since $\pi\circ\varphi\in \Rep(A,H)$, $\pi(D)E=\pi\circ\varphi(A)E$ 
is dense in $E$. 
\par\smnoind
(b) Let $D=\phi(A)\subseteq M(B)$ and $\varphi$ be the canonical map from
$A$ to $D$. 
By part (a), $D$ is generated by $(\id\otimes \varphi)(X)$. 
Now, consider $\pi$ to be a faithful non-degenerate representation of
$B$ on a Hilbert space $H$. 
Then $\pi$ induces a faithful non-degenerate representation $\tilde\pi$
of $M(B)$ and hence of $D$. 
Let $Y=(\id\otimes \pi\circ\phi)(X)=(\id\otimes \tilde\pi\circ\varphi)(X) 
\eta (C\otimes {\cal K}(H))$. 
Since both $\pi(B)$ and $\tilde\pi(D)$ are generated by $Y$, we have, 
by Lemma \ref{1.13}, $\pi(B) = \tilde\pi(D)$ and by the injectivity of 
$\tilde\pi$, $B=D$. 
\par\medskip
\par\medskip
\par\medskip
 
\section{The universal duality}
\par\bigskip

In this section, we will study the universal object corresponding to unitary corepresentations of a Hopf $C^*$-algebra (which is an analogue of the full group $C^*$-algebras).
Throughout this section, $(S,\delta)$ is a Hopf $C^*$-algebra. 
\par\medskip

\begin{defn}\label{2.1} 
Let $B$ be a $C^*$-algebra and $H$ be a Hilbert space. 
A unitary $v\in M(B\otimes S)$ is called a \emph{unitary corepresentation of $S$ in $B$} (respectively, \emph{unitary corepresentation of $S$ on $H$} when $B={\cal K}(H)$) if $(\id\otimes \delta)(v) = v_{12}v_{13}$. 
A unitary corepresentation $v$ of $S$ on $H$ is said to be 
\emph{cyclic} if the $\varphi_v(S^*)$ is a cyclic representation on $H$ 
(where $\varphi_v(f) = ({\rm id}\otimes f)(v)$ for $f\in S^*$).
\end{defn}
\par\medskip

We can define a sort of \emph{universal dual} of $S$ as a pair 
$(B,v)$ where $v$ is a unitary corepresentation of $S$ in $B$ such that for any 
unitary corepresentation $u$ of $S$ in any $C^*$-algebra $D$, 
there exists a unique non-degenerate *-homomorphism $\varphi$ from 
$B$ to $M(D)$ with $u=(\varphi\otimes \id)(v)$. 
However, this type of objects is too abstract to handle
and the following ``stronger version'' seems more useful
(see e.g. \cite{BS} or \cite{Ng1}).
\par\medskip

\begin{defn}\label{2.2} 
$(\hat S_p, V_S)$ is said to be a 
\emph{strong dual} of $S$ if $V_S$ is a unitary corepresentation 
of $S$ in $\hat S_p$ such that 
\par\noindent
(1) for any unitary corepresentation $u$ of $S$ on a Hilbert space
$H$, there exists a non-degenerate representation $\varphi_u$ of 
$\hat S_p$ on $H$ such that $u=(\varphi_u\otimes \id)(V_S)$;
\par\noindent
(2) the set $\{(\id\otimes f)(V_S): f\in S^*\}\cap \hat S_p$ is dense 
in $\hat S_p$. 
\end{defn}
\par\medskip

\begin{rem}\label{2.3} 
(a) If the strong dual of $S$ exists, it is unique up to isomorphism. 

\smnoind
(b) It is clear that $S\tilde{\ } = \{ f\in S^*: (\id\otimes f)(V_S) \in \hat S_p\}$ is a closed ideal of $S^*$. 
\end{rem}

\medskip

The next proposition characterizes the existence of strong duals
(c.f. \cite[2.1]{Ng2}). 
It requires the following well known lemma.

\medskip

\begin{lem} \label{1.2} Let $B$ be an algebra. 
Then there exists an injection from the collection of 
equivalent classes of cyclic representations of $B$ (recall that a representation $\phi$ of $B$ on a Hilbert space $H$ is \emph{cyclic} if $\phi(B)\xi$ is dense in $H$ for some $\xi\in H$) to the 
set of positive sesquilinear forms on $B$. 
\end{lem}

\medskip

\begin{prop}\label{2.5} 
The strong dual of $S$ exists if and
only if there exists an ideal $N$ of $S^*$ such that for any
unitary corepresentation $X$ of $S$ on $H$, the following conditions
are satisfied:
\par\noindent 
(i) $\hat S_X = \overline{\{ (\id\otimes f)(X): f\in N \}}$ 
is a non-degenerate $C^*$-subalgebra of ${\cal L}(H)$;
\par\noindent 
(ii) $X\in M(\hat S_X\otimes S)$. 
\end{prop}
\par
\noindent {\bf Proof:} Suppose that the strong dual of $S$ exists and
let $N=S\tilde{\ }$ (Remark \ref{2.3}(b)). 
Then for any unitary corepresentation $X$ of $S$ on $H$, there exists
a representation $\varphi$ of $\hat S_p$ on $H$ such that
$X=(\varphi\otimes \id)(V_S)$. 
Thus, we have $\hat S_X = \varphi(\hat S_p)$ and $X\in M(\hat S_X\otimes S)$. 
Conversely, suppose that there exists $N\subseteq S^*$ 
satisfying the two conditions.
By definition, cyclic unitary corepresentations of $S$ induce cyclic
representations of $S^*$. 
Therefore, Lemma \ref{1.2} ensures that the collection of all 
equivalent classes of cyclic unitary corepresentations is a set.
It is now valid to consider the direct sum $U$ of all
equivalent classes of cyclic unitary corepresentations of $S$. 
It is obvious that $U$ is a unitary corepresentation.  
For any given cyclic unitary corepresentation
$X$ of $S$ on $H$, there exists a representation $\varphi$ of 
$\hat S_U$ on $H$ such that $X = (\varphi\otimes \id)(U)$. 
On the other hand, for any unitary corepresentation $Y$ of $S$ on $H$, 
since $N$ is an ideal of $S^*$, we see that $\{(\id\otimes g)(Y): g\in S^*\}\subseteq M(\hat S_Y)$. 
Hence $Y$ is cyclic if and only if $\hat S_Y$ acts cyclicly on $H$.
Since $\hat S_Y$ is a $C^*$-algebra, the canonical representation of 
$\hat S_Y$ decomposes into a direct sum of cyclic representations of 
$\hat S_Y$. 
It is clear that each of these direct summands corresponds to a cyclic
unitary corepresentation of $S$.
Thus, there exists a representation $\varphi$ of $\hat S_U$ on $H$ 
such that $Y=(\varphi\otimes \id)(U)$.
\par\medskip

\begin{prop}\label{2.6} 
Suppose that the strong dual 
$\hat S_p$ of $S$ exists. 
There exists a full comultiplication $\hat\delta_p:\hat S_p\rightarrow M(\hat S_p\otimes_{\max}\hat S_p)$ such that $(\hat\delta_p\otimes \id)(V_S) = 
(V_S)_{13}(V_S)_{23}\in M((\hat S_p\otimes_{\max}\hat S_p)\otimes S)$.
Moreover, there exists a coidentity for this comultiplication. 
If, in addition, $S\tilde{\ }$ is $S$-invariant (see Remark \ref{2.3}(b)), then $\hat S_p$ is a full Hopf $C^*$-algebra. 
\end{prop}

\medskip

The proofs of the first two statements are easy and the proof of the last one is similar to that of Proposition \ref{3.3} in the next section.
Moreover, by a similar argument as that of Proposition \ref{3.3}, if the set $\{ (g\otimes \id)(V_S): g\in \hat S_p^* \}\cap S$ is dense in $S$, then the converse of the last statement in the Proposition \ref{2.6} also holds. 

\medskip

\begin{defn}\label{2.7} 
A Hopf $C^*$-algebra $S$ is said to be
\par\smnoind
(a) \emph{dualizable} if the strong dual $(\hat S_p, V_S)$ exists and 
$S\tilde{\ }$ is $S$-invariant;
\par\smnoind
(b) \emph{symmetrically dualizable} if it is dualizable and 
the set $\{ (g\otimes \id)(V_S): g\in \hat S_p^* \}\cap S$ is
dense in $S$. 
\end{defn}
\par\medskip

\begin{eg} \label{2.8}
(a) If $V$ is a $C^*$-multiplicative unitary on $H$ (\cite[2.3(b)]{Ng2}), 
then the Hopf $C^*$-algebra $S_V = \overline{\{(\omega\otimes \id)(W):\omega\in{\cal L}(H)_*\}}$ is symmetrically dualizable. 
This includes those Hopf $C^*$-algebras studied in \cite{Baaj}, \cite{BS} and \cite{Wor3}.
Moreover, as the ``locally compact quantum groups'' studied in \cite{KV1} and \cite{KV} are defined by manageable multiplicative unitaries, they are also symmetrically 
dualizable. 
In particular, $C_0(G)$ and $C_r^*(G)$ are symmetrically dualizable. 
On the other hand, $C^*(G)$ is also symmetrically dualizable and more generally,
$(S_V)_p$ is dualizable if $V$ is a regular multiplicative unitary (see \cite{BS}). 
\par\smnoind
(b) Let $S$ be a unital $C^*$-algebra with the trivial comultiplication
$\delta(a) = 1\otimes a$. 
Then $S$ is dualizable with $\hat S_p = \bf C$. 
$S$ is certainly not symmetrically dualizable unless $S=\bf C$. 
\par\smnoind
(c) Suppose that $S$ and $T$ are two Hopf $C^*$-algebras with 
coidentities $e_S$ and $e_T$ respectively. 
If both $S$ and $T$ are dualizable (respectively, symmetrically 
dualizable), then so is $S\otimes T$. 
In fact, if $Z$ is a unitary corepresentation of $S\otimes T$, then 
$X=(\id\otimes \id\otimes e_T)(Z)$ and $Y=(\id\otimes e_S\otimes \id)(Z)$ 
are unitary corepresentations of $S$ and $T$ respectively such that
$X_{12}Y_{13} = Y_{13}X_{12}$. 
\end{eg}
\par\medskip

\begin{rem}
Given a coaction $\epsilon$ of a Hopf $C^*$-algebra $S$ on a $C^*$-algebra $A$, one can define covariant representations for $\epsilon$.
Moreover, similar to the definition of the strong dual, one can also define the \emph{full crossed product} of $\epsilon$ (which may or may not exist). 
It can be shown that if $S$ is dualizable, then the full crossed product of $\epsilon$ exists if and only if there exists a covariant representation for $\epsilon$. 
However, as these are not relevant to the main concern of this paper, we will not give the details here. 
\end{rem}
\par\medskip
\par\medskip
\par\medskip

\section{Fourier duality and representations}
\par\bigskip

In this section, we will study a general duality framework. 
\par\medskip

\begin{defn}\label{3.1} 
Let $S$ and $T$ be Hopf $C^*$-algebras and $X\in M(T\otimes S)$ be a unitary. 
Then $(T,X,S)$ is called a \emph{Fourier duality} if 
\begin{equation}
\label{birep}
(\id\otimes \delta_S)(X) = X_{12}X_{13} \quad {\rm and} \quad 
(\delta_T\otimes \id)(X) = X_{13}X_{23}
\end{equation} 
as well as 
$$\overline{\{(\id\otimes f)(X): f\in S^*\}\cap T} = T
\quad {\rm and} \quad \overline{\{(g\otimes \id)(X): g\in T^*\}\cap S} =S.$$ 
\end{defn}
\par\medskip

\begin{eg}\label{3.2}
(a) For any $C^*$-multiplicative unitary $V$ (\cite[2.1]{Ng2}), we recall that the following algebras $\hat S_V = \overline{\{(\id\otimes \omega)(W):\omega\in{\cal L}(H)_*\}}$ and $S_V = \overline{\{(\omega\otimes \id)(W):\omega\in{\cal L}(H)_*\}}$ are Hopf $C^*$-algebras. 
Furthermore, $(\hat S_V, V, S_V)$ is a Fourier duality. 
In particular, it is the case for regular (\cite[3.3]{BS}), 
semi-regular and balanced (\cite[3.1]{Baaj}), or manageable (\cite[1.2]{Wor3} and also \cite{KV}) 
multiplicative unitaries. 

\smnoind
(b) If $(T, X, S)$ and $(T', X', S')$ are two Fourier dualities, 
then so is $(T\otimes T',(X\otimes X')^{\sigma_{23}},S\otimes S')$
where $\sigma_{23}$ is the flips of the two middle variables. 
In particular, if $V$ and $W$ are multiplicative unitaries such 
that $V$ is regular and $W$ is manageable, then $(\hat S_V\otimes 
\hat S_W, (V\otimes W)^{\sigma_{23}}, S_V\otimes S_W)$ is also a Fourier
duality. 
However, there is no reason to believe that $(V\otimes W)^{\sigma_{23}}$
is either regular or manageable. 

\smnoind
(c) If $S$ is symmetrically dualizable (see Definition \ref{2.7}(b)), then 
$(\hat S_p, V_S, S)$ is a Fourier duality. 

\smnoind
(d) Let $V$ and $W$ be two regular multiplicative unitaries. 
Then $((\hat S_V)_p\otimes \hat S_W, (V_{S_V}\otimes W)^{\sigma_{23}}, 
S_V\otimes S_W)$ is also a Fourier duality. 
In general, it is not possible to represent this sort of duality using multiplicative unitary unless $V$ is amenable. 

\smnoind
(e) Suppose that $(T,X,S)$ is a Fourier duality. 
Then $(S, X^\sigma, T)$, $(T^{op}, X^*, S^{op})$ and 
$(S^{op}, X^\top, T^{op})$ are all Fourier dualities 
(where $\sigma$ is the flip of variables, $X^\top = (X^\sigma)^*$ and 
$S^{op}$ is the Hopf $C^*$-algebra which has the same 
underlying $C^*$-algebra as $S$ but with the comultiplication 
$\delta_{S^{op}} = \sigma\circ\delta_S$). 
Note that in the setting of multiplicative unitaries, it is only 
valid to talk about $X$ and $X^\top$ and 
the other two only appear when irreducibility is considered. 
\end{eg}
\par\medskip

Suppose that $S$ and $T$ are pre-Hopf $C^*$-algebras (see Definition \ref{1.6}(b)) and $X\in M(T\otimes S)$ satisfying equalities (\ref{birep}).
Let $j_S$ and $j_T$ be the canonical homomorphisms from 
$T^*$ and $S^*$ to $M(S)$ and $M(T)$ respectively i.e. 
$j_S(g) = (g\otimes \id)(X)$ and $j_T(f) = (\id\otimes f)(X)$. 
We denote 
$$S^\# = \{ f\in S^*: j_T(f)\in T\}\qquad {\rm and}\qquad 
T^\# = \{ g\in T^*: j_S(g)\in S\}.$$
It is clear that $S^\#$ and $T^\#$ are closed ideals of $S^*$ 
and $T^*$ respectively. 
\par\medskip

\begin{prop}\label{3.3} 
Let $T$ and $S$ be pre-Hopf $C^*$-algebras. 
Suppose that $X\in M(T\otimes S)$ satisfies the conditions in Definition \ref{3.1}.
Then $S^\#$ is $S$-invariant if and only if $T$ is a Hopf
$C^*$-algebra. 
Consequently, for any Fourier duality $(T,X,S)$, $S^\#$ and $T^\#$ 
are automatically $S$-invariant and $T$-invariant respectively. 
\end{prop}
\par
\noindent {\bf Proof:}
Suppose that $T$ is a Hopf $C^*$-algebra. 
For any $f\in S^\#$ and $g\in T^\#$,
$$j_T(f\cdot j_S(g)) = (\id\otimes f\cdot j_S(g))(X) = 
(g\otimes \id \otimes f)(X_{13}X_{23}) = 
(g\otimes \id)(\delta_T(j_T(f)))$$
which is in $T$. 
Hence, $f\cdot j_S(g)\in S^\#$ and so $S^\#\cdot S \subseteq S^\#$ 
(by the density condition in Definition \ref{3.1}). 
Similarly, we have $j_T(j_S(g)\cdot f)=(\id\otimes g)(\delta_T(j_T(f)))$ and $S\cdot S^\# \subseteq S^\#$.
Conversely, suppose that $S^\#$ is $S$-invariant. 
For any $g'\in S^\#$, there exist $s\in S$ and $g\in S^\#$ 
such that $g'=s\cdot g$. 
Moreover, for any $t\in T$, $X(t\otimes s)\in T\otimes S$ can be approximated by elements of the form 
$\sum b_i\otimes a_i$ ($a_i\in S$; $b_i\in T$) and 
$$\delta_T((\id\otimes g')(X))(1\otimes t) =  (\id\otimes \id 
\otimes g)(X_{13}X_{23}(1\otimes t\otimes s))$$ which can be  
approximated by $\sum (\id\otimes b_i\cdot g)(X)\otimes a_i$.
This shows that $\delta_T((\id\otimes g')(X))(1\otimes t)\in
T\otimes T$ and thus, $\delta_T(T)(1\otimes T)\subseteq T\otimes T$.
On the other hand, by considering $g'=g\cdot s$, we see that $(t\otimes 1)\delta_T((\id\otimes g')(X))\in T\otimes T$ and so, $(T\otimes 1)\delta_T(T)\subseteq T\otimes T$. 
\par\medskip

\begin{cor}\label{3.4}
If $(T, X, S)$ is a Fourier duality, then $S^\#$ is invariant under the usual conjugation (which is given by $f^*(x) = \overline{f(x^*)}$). 
\end{cor}
\par
\noindent {\bf Proof:}
Since $S^\#$ is $S$-invariant, by \cite[III.4.4]{Tak}, any $f\in S^\#$ can be approximated by elements of the form 
$\sum_k s_k\cdot g_k$ where $s_k\in S$ and $g_k\in S^\#\cap S^*_+$.
Therefore, $f^*$ can be approximated by $\sum_k g_k\cdot s_k^*\in S^\#$. 
\par\bigskip

If, in addition, $S$ is dualizable, then we have two ``dualities'': 
$(\hat S_p, V_S, S)$ and $(T, X, S)$. 
It is natural to ask whether there is any relation 
between $T$ and $\hat S_p$.
By definition, there exists a $*$-homomorphism $q_T$ from 
$\hat S_p$ to $M(T)$ such that $X=(q_T\otimes \id)(V_S)$. 
It is not hard to show that $q_T$ is a Hopf \mbox{$*$-homomorphism}. 
Moreover, the following result shows that the image of $q_T$ is exactly $T$
(which is a bit surprising in this general situation of Fourier duality). 
\par\medskip

\begin{prop}\label{3.5}
Let $(T,X,S)$ be a Fourier duality. 
If $S$ is dualizable and $q_T$ is the induced map from 
$\hat S_p$ to $M(T)$, then $S$ is symmetrically dualizable (i.e. $(\hat S_p, V_S, S)$ is a Fourier duality) and $T = q_T(\hat S_p)$. 
\end{prop}
\par
\noindent {\bf Proof:}
$S$ is symmetrically dualizable because
$\{ (h\otimes \id)(V_S): h\in q_T^*(T^\#)\}$ is dense in $S$. 
Now, suppose that $A=q_T(\hat S_p)$.
Since both $S\tilde{\ }$ (see Remark \ref{2.3}(b)) and $S^\#$ are ideals of $S^*$, 
$S^\#\cdot S\tilde{\ }\subseteq S\tilde{\ }\cap S^\#$. 
Therefore, $T\cdot A \subseteq A\cap T$. 
Moreover, since $q_T$ is non-degenerate, $T\cdot A = T$. 
This implies that $T\subseteq A$. 
On the other hand, it is clear that $X=(q_T\otimes \id)(V_S)\in M(A\otimes S)$. 
Take any $h\in S^\#$ and any $s\in S$ such that $h(s)=1$. 
For any $a\in A$, because $X(A\otimes S) = A\otimes S$, we can approximate 
$a\otimes s$ by elements of the form 
$$\sum_i X((\id\otimes f_i)(X)
\otimes s_i) = \sum_i (\id\otimes \id\otimes f_i)(X_{12}X_{13})
(1\otimes s_i)$$ 
(where $f_i\in S\tilde{\ }$ and $s_i\in S$). 
Hence, $a$ can be approximated by elements of the form 
$$\sum_i (\id\otimes s_i\cdot h)(X)(\id\otimes f_i)(X)\in 
T\cdot A = T.$$
 
\par\bigskip

The next result gives a relation between 
Fourier duality and $C^*$-algebra generated by affiliated element.
It is a generalization of \cite[1.6.6]{Wor3}. 
\par\medskip

\begin{prop}\label{3.6}
Let $(T,X,S)$ be a Fourier duality such that both $S$ and $T$ are
separable. 
Then $S$ is generated by $X$ in the sense of Woronowicz 
(see Definition \ref{1.11}) and $T$ is generated by $X^\sigma$.
\end{prop}
\par
\noindent {\bf Proof:}
Let $H$ be a Hilbert space. 
Let $\varphi_1, \varphi_2 \in \Rep(S,H)$ such that $(\id\otimes \varphi_1)(X) = (\id\otimes \varphi_2)(X)$. 
Since $\{ f\in T^\#:(f\otimes \id)(X)\}$ is dense in $S$, we have
$\varphi_1 = \varphi_2$. 
This shows that $X$ satisfies condition I of \cite[4.2]{Wor2}. 
Next, suppose that $B\in C^*(H)$ and $\pi\in \Rep(S,H)$ such that $(\id\otimes \pi)(X)\eta (T\otimes B)$. 
Since $X$ is bounded, $(\id\otimes \pi)(X)\in M(T\otimes B)$. 
Moreover, as $X$ is a unitary, $(\id\otimes \pi)(X)(T\otimes B)$ is
dense in $T\otimes B$. 
We now take any $\omega\in T^\#\setminus (0)$.
As $T^\#$ is $T$-invariant, $T\cdot\omega\subseteq T^\#$ and hence
$(\omega\otimes \id)(X(t\otimes 1)) = (t\cdot\omega\otimes \id)(X)\in S$. 
By putting $F=T$ and $r=X$ in \cite[4.2(II)]{Wor2}, we see that $S$ is
generated by $X$. 
Similarly, $T$ is generated by $X^\sigma$
\par\bigskip

The argument in Proposition \ref{3.6} also shows that if $S$
is separable and dualizable, then 
$\hat S_p$ is generated by $V_S^\sigma$. 
Hence, if $S$ and $T$ are separable, we can also obtain 
Proposition \ref{3.5} by Propositions \ref{1.14}(b) and \ref{3.6}.
\par\medskip

\begin{rem}
Suppose that $(T,X,S)$ is a Fourier duality. 

\smnoind
(a) It may not be true in general that $S^\#$ separates points of $S$ 
but we can always ``quotient out the annihilators of $S^\#$''.
More precisely, if $I=\bigcap_{f\in S^\#}\ker f$, then it can be shown that $S_0 = S/I$ is a Hopf $C^*$-algebra and $(T, (\id\otimes q)(X), S_0)$ is a Fourier duality (where $q$ is the quotient map from $S$ to $S_0$) with $q^*(S_0^\#) = S^\#$. 
The same can be done on $T$. 

\smnoind
(b) It is clear that $j_T(S^\#)$ is a subalgebra of $T$. 
An interesting question is when it will be a $*$-subalgebra. 
Let us first consider the (unbounded) ``coinvolution'' 
$\kappa_S$ from $j_S(T^\#)$ to $S$ given by
$$\kappa_S(j_S(g)) = (g\otimes \id)(X^*) $$
$(g\in T^\#)$ (it is well defined because of Corollary \ref{3.4}). 
Suppose that $S^\#$ and $T^\#$ separate points of $T$ and $S$ respectively. 
Then $S$ can be regarded as a subalgebra of $(S^\#)^*$. 
Moreover, we have the following equivalent statements: 

\smnoind
(i) $j_T(S^\#)$ is a $*$-subalgebra of $T$.

\smnoind
(ii) $\kappa_S$ can be extended to a isometric anti-isomorphism from 
$S$ onto $S$.

\smnoind
(iii) $\kappa_S$ can be extended to a weak-*-continuous anti-isomorphism 
from $S^{\# *}$ onto itself.

\smnoind
Again, as this is not our main concern, we will not give the proof here. 
\end{rem}

\medskip

Next, we turn to representations of a Fourier duality $(T,X,S)$.
It is natural to start with unitary corepresentations of 
$S$ and $T$ and look for compatibility conditions with $X$. 
\par\medskip

\begin{defn}\label{4.1}
Let $(T, X, S)$ be a Fourier duality. 

\smnoind
(a) If $u$ and $v$ are unitary corepresentations of $T$ and $S$
respectively in $B$, then $(u,v)$ is said to be a \emph{unitary 
$X$-covariant representation} if it satisfies the following form of commutation relation:
$$u^\sigma_{12}X_{13}v_{23} = v_{23}u^\sigma_{12}.$$ 

\smnoind
(b) If $\mu$ and $\nu$ are representations of $S$ and $T$ respectively
on $H$, then $(\mu,\nu)$ 
is said to be a \emph{$X$-covariant
representation} if $((\mu\otimes \id)(X^\sigma),(\nu\otimes \id)(X))$ 
is a unitary $X$-covariant representation. 
Moreover, $(\mu,\nu)$ is said to be \emph{faithful} if and only if both 
$\mu$ and $\nu$ are injective. 
\end{defn}
\par\medskip

\begin{rem}\label{4.2} 
(a) Note that if $X$ is any unitary in $M(T\otimes S)$ such that there 
exist unitary corepresentations $u$ and $v$ of $T$ and $S$ 
respectively satisfying the equality in Definition \ref{4.1}(a), then
\cite[1.5]{Ng2} implies that $X$ satisfies equalities (\ref{birep}). 

\smnoind
(b) A unitary $X^\sigma$-covariant representation $(v, u)$ means that 
$v^\sigma_{12}X^\sigma_{13}u_{23} = u_{23}v^\sigma_{12}$ or
equivalently, $v_{23}X_{13}u^\sigma_{12} = u^\sigma_{12}v_{23}$.
In general, there is no relation between $X$-covariant
representations and \linebreak $X^\sigma$-covariant representations.

\smnoind
(c) Let $u$ and $v$ be unitary corepresentations of $S$ and $T$ 
respectively on a Hilbert space $H$. 
Then $(u,v)$ is a unitary $X$-covariant (respectively, $X^*$-covariant) 
representation if and only if $(v,u)$ is a unitary $X^\top$-covariant 
(respectively, $X^\sigma$-covariant) representation. 
\end{rem}
\par\medskip

\begin{prop} \label{4.3} 
If $(\mu,\nu)$ is a $X$-covariant representation, then 
$W=(\nu\otimes \mu)(X)$ is a multiplicative unitary. 
Moreover, both $\mu(S)$ and $\nu(T)$ are Hopf $C^*$-algebras (with coactions induced by $W$), both $\mu$ and $\nu$ are Hopf $*$-homomorphisms and $W\in M(\nu(T)\otimes \mu(S))$ satisfies equalities (\ref{birep}). 
\end{prop}
\par\medskip

In fact, since $(\mu\otimes \mu)\delta_S(s) = W(\mu(s)\otimes 1)W^*$, 
the map $\bar \delta_S$ on $\mu(S)$ defined by 
$\bar \delta_S(\mu(s)) = (\mu\otimes \mu)\delta_S(s)$
is a well defined comultiplication. 

\bigskip

Hence, $(\mu,\nu)$ can be regarded as a sort of multiplicative 
unitary realization of the abstract Fourier duality $(T, X, S)$. 
However, $\mu(S)$ and $\nu(T)$ may not be defined by $W$ in the usual
way i.e. it is possible that $\nu(T)\neq \hat S_W$
or $\mu(S)\neq S_W$ (see Example \ref{3.2}(a)). 

\bigskip

On the other hand, there may not exist any faithful $X$-covariant
representation (e.g. if $(C^*(G), X, C_0(G))$ is the strong dual of $S=C_0(G)$
where $G$ is a non-amenable locally compact group). 

\par\medskip
\par\medskip
\par\medskip

\section{Reduced duality and uniqueness}
\par\bigskip

In this section, we will define reduced Fourier dualities and study their uniqueness. 
This is done by using a notion of Fourier algebra of an arbitrary Hopf $C^*$-algebra.
For a ``nice'' multiplicative unitary $V$, we can recover $S_V$ from the structure of $(S_V)_p$ through this construction without looking at $V$ (this fact is important for later applications in \cite{Ng3}, \cite{Ng4} and \cite{Ng5}).

\bigskip 

For simplicity, we call a subspace $N$ of $S^*$ a 
\emph{Hopf ideal} (respectively, \emph{left Hopf ideal} or \emph{right Hopf
ideal}) if it is a non-trivial $S$-invariant closed 2-sided ideal (respectively,
left ideal or right ideal) of $S^*$. 
For any Fourier duality $(T,X,S)$, the subspace $S^\#$ is a Hopf ideal. 
Note that we use the term ``Hopf ideal'' in a different way from the literatures on Hopf algebras (although the quotient of $S$ with respect to the annihilators of
an Hopf ideal is a Hopf $C^*$-algebra). 
\par\medskip

\begin{defn}\label{5.1}
(a) The \emph{Fourier algebra}, $A_S$, of a Hopf $C^*$-algebra $S$ is the 
intersection of all Hopf ideals of $S^*$. 

\smnoind
(b) A Hopf $C^*$-algebra $S$ is said to be \emph{irreducible} if any Hopf ideal of $S^*$ separates points of $S$. 

\smnoind
(c) A Fourier duality $(T,X,S)$ is said to be \emph{reduced} 
if $\{ (\id\otimes f)(X): f\in A_S\}$ is dense in $T$ and 
$\{ (g\otimes \id)(X): g\in A_T\}$ is dense in $S$. 
\end{defn}
\par\medskip

Note that the definition of the Fourier algebras comes from the corresponding property of $A(G)$ (see \cite[2.5]{DR}).
In the following, we will show that it is really a ``right definition''. 

\medskip

\begin{lem}\label{5.2}
Let $S$ be a irreducible Hopf $C^*$-algebra. 
Then its comultiplication $\delta_S$ is injective. 
\end{lem}
\par
\noindent {\bf Proof:}
Let $I$ be the kernel of $\delta_S$. 
It is clear that $I^\bot = \{ f\in S^*: f(I) = 0\}$ is a 
$S$-invariant closed subspace of $S^*$. 
Moreover, for any $f,g\in S^*$, we have $fg\in I^\bot$ and in 
particular, $I^\bot$ is an ideal of $S^*$. 
Hence, either $I^\bot = (0)$ (i.e. $I=S$ and is impossible as $\delta$ is a non-degenerate $*$-homomorphism) or $I^\bot$ is 
$\sigma(S^*,S)$-dense in $S^*$ which implies that $I^\bot =S^*$ (i.e. $I=(0)$).
\par\bigskip

It is natural to ask whether we can get a sort of ``irreducible quotient'' of a 
general Hopf $C^*$-algebra. 
If $S$ is a Hopf $C^*$-algebra such that $A_S \neq (0)$, we can take the quotient
$S_{irr}=S/(A_S)^\top$ with the canonical quotient map $Q_{S}$.
\par\medskip

\begin{prop}\label{5.3}
Let $S$ and $T$ be Hopf $C^*$-algebras such that $A_S\neq (0)$ and $A_T\neq (0)$.
Then $S_{irr}$ is an irreducible Hopf $C^*$-algebra with 
$A_{S_{irr}}= A_S$. 
Moreover, if $(T,X,S)$ is a Fourier duality and 
$X_{irr}=(Q_{T}\otimes Q_{S})(X)$, then 
$(T_{irr}, X_{irr}, S_{irr})$ is a reduced Fourier duality. 
Consequently, if both $S$ and $T$ are irreducible with non-zero Fourier algebras, then $(T,X,S)$ is automatically reduced.
\end{prop}
\par
\noindent {\bf Proof:}
By Remark \ref{1.8}(b) (and the fact that Fourier algebras are invariant and hence are defined by central projections), $S_{irr}$ and $T_{irr}$ are both Hopf $C^*$-algebras. 
Since $S_{irr}^* \cong \overline{A_S}^{\sigma(S^*,S)}$, 
we have $A_{S_{irr}} = A_S$. 
It is clear that $Q_{S}$ and $Q_{T}$ are Hopf 
$*$-homomorphisms and $X_{irr}$ satisfies equalities (\ref{birep}). 
It is required to show that $(T_{irr}, X_{irr}, S_{irr})$ is a 
reduced Fourier duality. 
Let 
$$N=\{ g\in T^*: (g\otimes f)(X)=0 \;\:{\rm for\;\: any}\;\: 
f\in A_S\}.$$
Then $N$ is a closed subspace of $T^*$. 
Moreover, for any $g\in N$, $f\in A_S$ and $h\in S^\#$, 
$$((j_T(h)\cdot g)\otimes f)(X)= (g\otimes fh)(X)=0$$ 
(since $A_S$ is an ideal of $S^*$) and similarly 
$(g\cdot j_T(h)\otimes f)(X) = (g\otimes hf)(X) = 0$. 
This shows that $N$ is $T$-invariant. 
On the other hand, for any $g\in N$, $l\in T^*$ and $f\in A_S$, 
$$(gl\otimes f)(X) = (g\otimes j_S(l)\cdot f)(X) = 0$$ 
(as $A_S$ is $S$-invariant and hence is $M(S)$-invariant) and 
$(lg\otimes f)(X)=0$. 
Thus, either $N=(0)$ or $N$ is a Hopf ideal (or equivalently, $A_T\subseteq N$). 
Suppose that $A_T\subseteq N$. 
Then $(g'\otimes f')(X_{irr})=0$ for all $g'\in A_T = A_{T_{irr}}$ and 
$f'\in A_S = A_{S_{irr}}$. 
This implies that $X_{irr}=0$ (as $A_{S_{irr}}$ and $A_{T_{irr}}$ 
separate points of $S_{irr}$ and $T_{irr}$ respectively and are invariant) and is 
impossible because $X_{irr}$ is a unitary. 
Therefore, $N=(0)$ and $\{ (\id\otimes f)(X): f\in A_S\}$ is norm dense 
in $T$. 
Similarly, $\{ (g\otimes \id)(X): g\in A_T\}$ is dense in $S$. 
\par\bigskip

From the proof of the above proposition, we see that if we consider 
the duality between $S^*$ and $T^*$ defined by $X$, then the ``polar'' of an
ideal of $T^*$ is $S$-invariant and the ``polar'' of a closed $T$-invariant
subspace is an ideal of $S^*$. 
\par\medskip

\begin{eg} \label{5.4}
(a) For any locally compact group $G$, $C^*_r(G)$ is an irreducible 
Hopf $C^*$-algebra (by \cite[2.4]{DR}). 
In fact, we have $C^*_r(G) = C^*(G)_{irr}$. 

\smnoind
(b) If $S=S_V$ (see Example \ref{3.2}(a)) for a regular irreducible multiplicative unitary $V$, then $(\hat S_p)_{irr} = \hat S_V$
(see Corollary \ref{5.10}(c) below). 
It is an interesting fact since we can obtain the ``reduced dual 
object'' $\hat S_V$ as well as
$V\in M(\hat S_V\otimes S_V)$ without knowing \emph{a priori} that 
$S$ comes from a Kac system nor assuming any special representation 
of $S$. 
\end{eg}
\par\medskip

In the following, we will give some uniqueness type results for
reduced Fourier dualities (Theorem \ref{5.13} and Corollary \ref{5.14}(b)). 
We will also consider their relations with multiplicative unitaries 
(Proposition \ref{5.12}(b) and Corollary \ref{5.14}(a)). 
\par\medskip

\begin{lem}\label{5.5}
Let $(\mu ,\nu)$ be a $X$-covariant representation and $e_\mu$ be the 
support projection of the extension $\tilde\mu$ of $\mu$ on $S^{**}$.
Then $\delta_S^{**}(e_\mu) \geq e_\mu\otimes 1$.
\end{lem}
\par
\noindent {\bf Proof:}
Since $(\mu,\nu)$ is a $X$-covariant representation, 
$$(\mu\otimes \id)\delta_S(x) = (\nu\otimes \id)(X)(\mu(x)\otimes 1)
(\nu\otimes \id)(X)^*$$ 
($x\in S$). 
Note that by the weak density of $S$ in $S^{**}$, the same equation 
holds for $x\in S^{**}$ (with $\mu$ replaced by $\tilde\mu$ and $\delta_S$ by $\delta_S^{**}$). 
Therefore, $(\tilde\mu\otimes \id)\delta_S^{**}(e_\mu) = 1$ and so 
$(\tilde\mu\otimes \id)((e_\mu\otimes 1)\delta_S^{**}(e_\mu)) = 
(\tilde\mu\otimes \id)(e_\mu\otimes 1)$. 
Since $\tilde\mu$ is injective on $e_\mu S^{**}$, we have 
$(e_\mu\otimes 1)\delta_S^{**}(e_\mu) = e_\mu\otimes 1$. 
\par\medskip

\begin{prop}\label{5.6}
Let $(T, X, S)$ be a Fourier duality and $p$ be a non-zero
central projection in $S^{**}$. 
Let $(\mu,\nu)$ be a $X$-covariant representation. 
If $\delta_S^{**}(p)\geq 1\otimes p$, then $p\geq e_\mu$. 
\end{prop}
\par
\noindent {\bf Proof:}
Let $W=(\nu\otimes \mu)(X)\in\nu(T)''\bar \otimes \mu(S)''$. 
Then, by Proposition \ref{4.3}, $(\tilde\mu\otimes \tilde\mu)
\delta_S^{**}(s) = W(\tilde\mu(s)\otimes 1)W^*$ for any $s\in S^{**}$ (where $\tilde \mu$ is as in the previous lemma). 
If $\tilde p=\tilde\mu(p)$, then $\tilde p$ is a central 
projection in $\mu(S)''$ and we have 
\begin{eqnarray*}
W((1-\tilde p)\otimes \tilde p)W^* & = &
W(1\otimes \tilde p)W^* - (1\otimes \tilde p)W(\tilde p\otimes 1)W^*\\ 
& = & 1\otimes \tilde p - (\tilde\mu\otimes \tilde\mu)((1\otimes p)
\delta_S^{**} (p))\quad = \quad 0.
\end{eqnarray*}
Hence either $\tilde p=0$ or $\tilde p=1$. 
Suppose that $\tilde\mu(p)=\tilde p=0$. 
Then $p\leq 1-e_\mu$ and, by Lemma \ref{5.5}, 
$1\otimes p + e_\mu\otimes 1 \leq \delta_S^{**} (p) + 
\delta_S^{**} (e_\mu) \leq 1\otimes 1$. 
Therefore, $e_\mu\otimes 1 \leq 1\otimes(1-p)$ and so 
either $p=0$ or $e_\mu=0$ which contradicts the 
assumption. 
Thus, we have $\tilde p=1$ or, equivalently, $p\geq e_\mu$. 
\par\medskip

\begin{rem}\label{5.7}
(a) The idea of the proof of Proposition \ref{5.6} is mainly from that of 
\cite[2.4]{DR} but since we consider central projections only, we have 
a simpler proof here. 

\smnoind
(b) We also have a corresponding result for $\nu$ in Lemma \ref{5.5} 
and Proposition \ref{5.6} (namely, $\delta_T^{**}(e_\nu)\geq 1\otimes e_\nu$ 
and $\delta_T^{**}(q)\geq q\otimes 1$ implies that $q\geq e_\nu$).
\end{rem}
\medskip

\begin{cor}\label{5.8}
Let $(T,X,S)$ be a Fourier duality such that $A_S\neq (0)$ and 
let $(\mu,\nu)$ be a $X$-covariant representation.
Then $e_\mu\cdot S^* \subseteq A_S$ and $\mu$ induces a representation 
$\mu_{irr}$ of $S_{irr}$ such that $(\mu_{irr},\nu)$ is a 
$(\id\otimes Q_{S})(X)$-covariant representation. 
\end{cor}
\par
\noindent {\bf Proof:}
Let $p$ be the central projection of $S^{**}$ such that $A_S = p\cdot S^*$ (as $A_S$ is $S$-invariant). 
Then by Propositions \ref{5.6} and \ref{1.9} (note that $A_S$ is a left ideal of $S^*$), $p\geq e_\mu$. 
Hence $e_\mu \cdot S^*\subseteq A_S$ and $\mu$ induces a representation 
$\mu_{irr}$ of $S_{irr}$.
It is easy to check that $(\mu_{irr},\nu)$ is a 
$(\id\otimes Q_{S})(X)$-covariant representation.
\par\bigskip

Therefore, if both $A_S$ and $A_T$ are non-zero, then any $X$-covariant
representation $(\mu, \nu)$ will induce a \mbox{$X_{irr}$-covariant}
representation. 
\par\medskip

\begin{thm} \label{5.9}
Let $(T,X,S)$ and $(S,Y,R)$ be Fourier dualities. 
Let $(\mu,\nu)$ and $(\rho,\lambda)$ be a $X$-covariant and 
a $Y$-covariant representations respectively.
Then 
$$e_\mu =e_{\lambda} \qquad {\rm as \ well \ as} \qquad A_S=e_\mu\cdot S^*=
\mu^*({\cal L}(H)_*)$$
and hence $A_S\neq (0)$. 
Moreover, $e_\mu\cdot S^*$ is the smallest left Hopf ideal as well as the
smallest right Hopf ideal of $S^*$. 
\end{thm}
\par
\noindent {\bf Proof:}
The first equality follows from Lemma \ref{5.5}, Proposition \ref{5.6} and Remark \ref{5.7}(b). 
Now, $e_\mu =e_{\lambda}$ implies that $e_\mu \cdot S^*$ is a Hopf ideal of $S^*$ 
(by Proposition \ref{1.9} and Lemma \ref{5.5}) and equals $A_S$ by Lemma \ref{5.8}. 
Moreover, $e_\mu \cdot S^* = e_\lambda \cdot S^*$ is contained in any left or right Hopf ideal of $S^*$ 
(by Propositions \ref{5.6} and \ref{1.9}). 
This completes the proof.
\par\medskip

\begin{cor}\label{5.10} 
(a) With the notations as in Theorem \ref{5.9}, $\mu$ (respectively, $\lambda$) induces a Hopf $*$-isomorphism from $S_{irr}$ to $\mu(S)$ (respectively, $\lambda(S)$). 

\smnoind
(b) Let $V$ and $W$ be two $C^*$-multiplicative unitaries such that 
$\hat S_W=S_V$ (see Example \ref{3.2}(a)). 
Then $S_V$ is irreducible with non-zero Fourier algebra. 

\smnoind
(c) Let $V$ be any regular irreducible multiplicative unitary and $(\hat S_p, W, S_V)$ be the strong dual of $S_V$.
If $(\mu,\nu)$ is any $W$-covariant representation,
then $\nu((\hat S_V)_p) \cong \hat S_V$ (as Hopf $C^*$-algebras). 
In particular, this applies to $(C^*(G), W_G, C_0(G))$. 
\end{cor}
\par
\noindent{\bf Proof:}
(a) Note that $\mu$ is a Hopf $*$-homomorphism from $S$ to $\mu(S)$ (Proposition \ref{4.3}).
Moreover, by Theorem \ref{5.9}, $\ker \mu = (\mu^*({\cal L}(H)_*)^\top = (A_S)^\top$. 
Hence, $\mu$ induces a Hopf $*$-isomorphism from $S_{irr}$ to $\mu(S)$.

\smallskip\noindent
(b) This part clearly follows from the fact that $A_{S_V}=
L_V^*({\cal L}(H)_*)$ (Theorem \ref{5.9}) which separates points of $S_V$. 

\smallskip\noindent
(c) Since $V$ is regular, $(\hat S_p, W, S_V)$ is a Fourier duality (by Proposition \ref{3.5}). 
Part (c) is now a direct application of parts (a) \& (b) (note that if $\hat q$ is the quotient from $\hat S_p$ to $\hat S_V$, then $(\hat q, \id)$ is also a $W$-covariant representation). 
\par\medskip

\begin{defn} \label{5.11} 
Suppose that $(T,X,S)$ is a Fourier duality. 
Then $(T,X,S)$ is said to be a \emph{Kac-Fourier duality} if both $S$ and $T$ are irreducible and there exist both a $X$-covariant representation and a $X^\sigma$-covariant
representation. 
\end{defn}
\par\medskip

\begin{prop}\label{5.12} 
Let $(T,X,S)$ be a Fourier duality. 

\smnoind
(a) If $(T,X,S)$ is a Kac-Fourier duality, then it is defined by a 
multiplicative unitary in the usual way. 

\smnoind
(b) $(T,X,S)$ is a Kac-Fourier duality if and only if there exist 
a $X^\sigma$-covariant representation as well as a faithful
$X$-covariant representation. 
It is the case if and only if $(T,X,S)$ and $(S,X^\sigma,T)$ are defined (in the usual 
sense) by two multiplicative unitaries $V$ and $W$ respectively.
\end{prop}
\par
\noindent 
{\bf Proof:}
(a) Suppose that $(\mu,\nu)$ is a $X$-covariant representation. 
Let $V$ be the multiplicative unitary $(\nu\otimes \mu)(X)$ (see Proposition \ref{4.3}). 
Then because $\mu^*({\cal L}(H)_*) = A_S$ and $\nu^*({\cal L}(H)_*) =A_T$ 
(by Theorem \ref{5.9}), we have, by Proposition \ref{5.3},
$$\mu(S) = \{ \mu((f\otimes \id)(X)): f\in A_T \} = \{ (\omega \circ \nu \otimes \mu)(X): \omega\in {\cal L}(H)_*\} = S_V$$ 
and $\nu(T) = \hat S_V$.  
Moreover, Corollary \ref{5.10}(a) implies that both $\mu$ and $\nu$ 
are injective. 
Therefore, $\mu$ and $\nu$ are Hopf $*$-isomorphisms from $S$ and $T$ to $S_V$ and $\hat S_V$ respectively (see Proposition \ref{4.3}). 

\smnoind
(b) As above, any $X$-covariant representation
of a Kac-Fourier duality is faithful. 
Conversely, if there exist a $X$-covariant representation as well as 
a $X^\sigma$-covariant representation such that one of them is 
faithful, then Corollary \ref{5.10}(a) shows that both $S$ and $T$ are
irreducible. 
Finally, one direction of the second statement clearly follows from part (a) while the other follows from the argument of Corollary \ref{5.10}(b).
\par\bigskip

Note that the multiplicative unitaries obtained above need not be 
semi-regular nor manageable. 
Hence Kac-Fourier duality is a far more general framework than Kac 
system. 

\bigskip

Finally, if $S$ is dualizable, we also have the following 
uniqueness type theorem for the ``dual object''. 
\par\medskip

\begin{thm} \label{5.13}
Let $(T,X,S)$, $(S,Y,R)$, $(\mu, \nu)$ and $(\rho, \lambda)$ 
be the same as in Theorem \ref{5.9}. 
Suppose that $S$ is dualizable. 
Then $\rho\circ q_R$ induces a Hopf $*$-isomorphism from $(\hat S_p)_{irr}$ to 
$\rho(R)$ (where $q_R$ is the canonical quotient map from $\hat S_p$ to $R$; see Proposition \ref{3.5}).
The same is true for $\nu\circ q_T$. 
Therefore, if both $\rho$ and $\nu$ are injective, then $R\cong T$ (as Hopf $C^*$-algebras) and $Y$ is the image of $X^\sigma$ under this isomorphism. 
\end{thm}
\par
\noindent {\bf Proof:}
Note that $(\mu, \nu\circ q_T)$ is a $V_S$-covariant representation 
while $(\rho\circ q_R, \lambda)$ is a $V_S^\sigma$-covariant 
representation (the strong dual $(\hat S_p, V_S, S)$ is a Fourier duality because of Proposition \ref{3.5}). 
Corollary \ref{5.10}(a) will now give the first two statements.
The final statement follows from Proposition \ref{3.5}.
\par\bigskip

The above theorem applies, in particular, to Kac-Fourier dualities. 
This also gives the following corollary. 
\par\medskip

\begin{cor}\label{5.14}
(a) Let $V$ and $W$ be two $C^*$-multiplicative unitaries 
(in particular, regular or manageable). 
If $S_V \cong \hat S_W$ as Hopf $C^*$-algebra, then 
$\hat S_V \cong S_W$ as Hopf $C^*$-algebras and $V$ is the image of $W^\sigma$ under these isomorphisms. 
Moreover, both $S_V$ and $\hat S_V$ are irreducible. 
In this case, if $U$ is another $C^*$-multiplicative 
unitary such that $S_V \cong S_{U}$, then $\hat S_V \cong 
\hat S_{U}$ and $U$ is the image of $V$ under these isomorphisms. 

\smnoind
(b) Let $V$ and $W$ be two semi-regular and balanced (or regular) 
multiplicative unitaries such that $V$ is irreducible.
If $S_V\cong S_W$ as Hopf $C^*$-algebras, then $\hat S_V\cong \hat S_W$
as Hopf $C^*$-algebras and $V$ is the image of $W$ under these isomorphisms. 
\end{cor}
\par\medskip
\par\medskip
\par

\par\medskip
\par\medskip
\par\medskip

\noindent Department of Pure Mathematics, Queen's University Belfast, Belfast BT7 1NN, United Kingdom.
\par\medskip
\noindent $e$-mail address: c.k.ng@qub.ac.uk

\end{document}